\documentclass[a4paper,11pt]{article}
\usepackage{amssymb}
\usepackage{amsfonts}
\usepackage{amsmath,tabu}
\usepackage{amsthm}
\usepackage{cite}
\usepackage{amsmath,tabu, color}

\newcommand{\R}{\mathbb{R}}

\newcommand{\N}{\mathbb{N}}

\newcommand{\cuad}{{\sqcap\kern-.68em\sqcup}}

\newcommand{\norm}[1]{\|#1\|}

\numberwithin{equation}{section}
\newtheorem{theorem}{Theorem}[section]

\newtheorem{proposition}[theorem]{Proposition}

\newtheorem{lemma}[theorem]{Lemma}
\newtheorem{corollary}[theorem]{Corollary}
\newtheorem{remark}[theorem]{Remark}
\newcommand{\bremark}{\begin{remark} \em}
\newcommand{\eremark}{\end{remark} }

\newcommand{\cH}{{\mathbb H}}

\newcommand{\cL}{{\mathcal L}}

\newcommand{\cN}{{\mathcal N}}

\allowdisplaybreaks[4]

\newcommand{\loglap}{L_{\text{\tiny $\Delta \,$}}\!}

\begin{document}
\begin{center}{\bf  \large Bounds of  Dirichlet eigenvalues for   Hardy-Leray operator }\medskip

 \bigskip

{\small

 {\sc  Huyuan Chen\footnote{chenhuyuan@yeah.net}}

\medskip
Department of Mathematics, Jiangxi Normal University,\\
Nanchang, Jiangxi 330022, PR China \\[10pt]

 {\sc  Feng Zhou\footnote{fzhou@math.ecnu.edu.cn} }

\medskip
Center for PDEs, School of Mathematical Sciences, East China Normal University, and NYU-ECNU Institute of Mathematical Sciences at NYU Shanghai,\\
Shanghai Key Laboratory of PMMP, Shanghai 200062, PR China \\[20pt]
 }

\begin{abstract}
The purpose of this paper is to study  the eigenvalues $\{\lambda_{\mu,i} \}_i$ for the Dirichlet Hardy-Leray operator, i.e.
$$ -\Delta u+\mu|x|^{-2}u=\lambda  u\ \  {\rm in}\ \,  \Omega,\quad\quad  u=0\ \ {\rm on}\  \ \partial\Omega,$$
where  $-\Delta  +\frac{\mu}{|x|^2}$ is the Hardy-Leray operator with $\mu\geq -\frac{(N-2)^2}{4}$ and $\Omega$ is a smooth bounded domain with $0\in\Omega$.
We provide lower bounds of $\{\lambda_{\mu,i} \}_i$, as well as   the Li-Yau's one when $\mu>-\frac{(N-2)^2}{4}$ and   Karachalios's bounds for $\mu\in [-\frac{(N-2)^2}{4},0)$. Secondly, we obtain Cheng-Yang's type upper bounds for $\lambda_{\mu,k}$. Addtionally we get the
 Weyl's limit of eigenvalues which is independent of the potential's parameter $\mu$. This interesting phenomena indicates that the inverse-square potential does not play an essential role for the asymptotic behavior of the spectrum  of the problem under study.

 \end{abstract}

\end{center}
 \vspace{1mm}


\vspace{1mm}

\setcounter{equation}{0}
\section{Introduction and main results}

Let $\Omega$ be  a smooth bounded domain in $\R^N$ with $N\ge 2$, $0\in\Omega$ and
 $\mu_0:=-\frac{(N-2)^2}{4}$ the best constant for the standard Hardy inequality.   The purpose of this paper is to study
 the bounds of  eigenvalues for the Dirichlet  problem
\begin{equation}\label{eq 1.1}
\left\{ \arraycolsep=1pt
\begin{array}{lll}
 \cL_\mu u=\lambda  u\quad \  &{\rm in}\ \   \Omega,\\[2mm]
 \phantom{ \loglap   }
  u=0\quad \ &{\rm{on}}  \ \, \partial\Omega,
\end{array}
\right.
\end{equation}
where  $\cL_\mu : = -\Delta  +\frac{\mu}{|x|^2}$ for $\mu\geq \mu_0$ is the Hardy-Leray operator.

When $\mu=0$, $\cL_\mu$ reduces to the Laplacian.  In 1912, Weyl   \cite{W} showed that  the $k$-th eigenvalue $\lambda_k $ of the Dirichlet problem for any smooth bounded domain $\Omega$ with the Laplacian operator has the following asymptotic behavior
\begin{equation}\label{W const}
\lim_{k\to+\infty}k^{-\frac2N}  \lambda_k= c_N |\Omega|^{-\frac2N} \quad {\rm with}\ \  c_N=(2\pi)^2|B_1|^{-\frac{2}{N}},
\end{equation}
here $c_N$ is named as Weyl's constant, $|\Omega|$ is the volume of $\Omega$ and $B_1$ is the unit ball centered at the origin.
 Later, P\'olya \cite{P} (in 1960) proved that
\begin{equation}\label{p-conj}
\lambda_{k} \geq C |\Omega|^{-\frac2N}k^{\frac2N}
\end{equation}
with $C=c_N, \, (N=2)$ and for any ``plane-covering domain"  in $\R^2$, (\,his proof also works in dimension $N\geq 3$)
 and he also conjectured that (\ref{p-conj}) holds with $C=c_N$ for any bounded domain in $\R^N$. Later on, Lieb \cite{L} proved (\ref{p-conj}) with some positive constant $C$ in a general bounded domain and then  Li-Yau \cite{LY} improved this constant in (\ref{p-conj}) to
\begin{equation}\label{LY const}
 C_N:=\frac{N }{N+2} c_N=\frac{N }{N+2}(2\pi)^{2 }  |B_1|^{-\frac{2}{N}}.
\end{equation}
The estimate (\ref{p-conj}) with $C=C_N$ is also called Berezin-Li-Yau's inequality because the constant $C_N$ is achieved thanks to Legendre transform in an earlier result obtained by Berezin \cite{Be}. Here we denote $C_N$ the Li-Yau's constant.  The Berezin-Li-Yau inequality then is generalized in \cite{CgW,L,K,M,CgY}   for degenerate elliptic operators.
  Upper bounds for the first $k$-eigenvalues  obtained in \cite{Ko} is controlled  by $c_N k^{\frac{2}{N}}$ together with lower order terms; later on, Cheng-Yang in \cite{CgY,CgY1,CgY2} developed  a very interesting  upper bound condition
  $$\lambda_{k}\leq \Big(1+\frac4N\Big)k^{\frac2N}\lambda_{1}, $$
 which  is also named Cheng-Yang's inequality.  More results on minimizing the eigenvalues problems could be found in  \cite{Fa,K,AE}, subject to the homogeneous Dirichlet boundary condition,  to \cite{A,BH} with  Neumann boundary condition.
  Our motivations of the Dirichlet  eigenvalues for Hardy operators are twofold, the one is the  Polya's conjecture  associating the zero of the Riemann Zeta function with the eigenvalue of a Hermitian operator and the latter is the role of the critical potential in the analysis of PDEs.

When $\mu\geq\mu_0$ and $\mu\not=0$,  the Hardy potential has homogeneity $-2$, which is critical from both
mathematical and physical viewpoints.   The Hardy-Leary problems  have numerous  applications in various
fields as molecular physics \cite{LL}, quantum cosmological models such as the
Wheeler-de-Witt equation (see e.g. \cite{BE}) and combustion models \cite{G}.
The Hardy inequalities
play a fundamental role in the study of Hardy-Leray problems. When the potential's  singularity  $\{0\}$ is in  $\Omega$,
the Hardy inequality \cite[(2)]{DD} (also see \cite{BrVa,BM}) reads as
\begin{equation}\label{Hardy inequalities}
 \int_\Omega |\nabla u(x) |^2 dx+ \mu_0\int_\Omega \frac{u^2(x)}{|x|^2}dx\geq c_0 \int_\Omega  u^2(x) dx,\quad \forall \, u\in H^1_0(\Omega)
 \end{equation}
for some constant $c_0>0$.  More Hardy-Lieb-Thirring inequalities inequalities could see \cite{EF,F,FL}.
The operator $\cL_\mu$  is then positive definite for $\mu\in [\mu_0,+\infty)$. More properties of Hardy-Leray operator
can be found in \cite{Caz2,MT2}.  It is worth noting that the inverse square potential plays an essential role
in isolated singularity of   elliptic Hardy equations.  Indeed, for a given `source'   $f$ defined in $\Omega$, the solutions of the Dirichlet problem
$$ \cL_\mu u=f \quad\text{in }\ \Omega\setminus\{0\},\qquad  u=0 \quad\text{on } \;\partial\Omega
 $$
with isolated singularities at $\{ 0 \}$ are classified fully  in \cite{CQZ}, thanks to a new formulation of distributional identity associated to some specific weight.  Extensive treatments of the associated semilinear problems are developed in \cite{AC,ChVe, ChVe1} via an introduction of a notion of very weak solution.

Due to the inverse square potential,   the spectral of the Hardy-Leray operator remains widely unexplored.    In \cite{D} the authors investigated the essential spectral for general Hardy-Leray potential,  for problem (\ref{eq 1.1}),  an increasing sequence of eigenvalues $\{\lambda_{\mu,i}(\Omega) \}_{i\in \N}$, simply denoted by $\{\lambda_{\mu,i} \}_{i\in \N}$ in the sequel, are obtained in \cite{VZ};   \cite{C} studied  the optimal first eigenvalue with respect to the domain,  and the lower bounds estimates of  eigenvalues for (\ref{eq 1.1}) in \cite{K1,K2} are derived by using the estimates of the related heat kernel and Sobolev inequalities. Precisely, the lower bounds state as following:
\begin{enumerate}
\item[(i)]  when $N\geq 3$ and $\mu_0< \mu<0$,
there holds that for $k\in\N$
\begin{equation}\label{test 1-0}
\lambda_{\mu,k}  \geq  \Big(1-\frac{\mu}{\mu_0}\Big)  \frac{N(N-2)}{4e}\omega_{_{N-1}}^{\frac2N}|\Omega|^{-\frac2N} k^{\frac{2}{N}};
\end{equation}
\item[(ii)] when $N\geq 3$ and $\mu=\mu_0$,
there holds that for $k\in\N$,
 $$ \lambda_{\mu_0,k}  \geq  e^{-1}S_N(N-2)^{-\frac{2(N-1)}{N}} \|\mathbb{X}_1^{\frac{1- N}2}\|_{L^2(\Omega)}^{-\frac4N} k^{\frac{2}{N}},$$
where $\mathbb{X}_1(x):=\big(-\log(\frac{|x|}{D_0})\big)^{-1},\ \forall\, x\in\Omega$,
$D_0:=\max_{x\in\partial\Omega}|x|$, and $S_N:=2^{2/N}\pi^{1+1/N}\Gamma(\frac{N+1}{2})$
is the best constant of Sobolev inequality in $\R^N$ and $\Gamma$ is the well-known Gamma function. \end{enumerate}
We note that the defect of the lower bound  (\ref{test 1-0})  is the coefficient $(1-\frac{\mu}{\mu_0})\frac{N(N-2)}{4e} \to 0$ as $\mu\to\mu_0^+$.\smallskip

It is natural to ask whether the Hardy-Leray potential plays an essential role in the estimates of the spectral estimates or how it works on
the related eigenvalues.  Our main aim in this paper is to answer this question and to establish the lower   and upper bounds,
and show the limit the eigenvalues as $k\to+\infty$. Now we state our lower bounds  as follows.

\begin{theorem}\label{teo 1.2}
Let  $D_0=\max_{x\in\partial\Omega}|x|$ and $\{\lambda_{\mu,i} \}_{i\in\N}$  be the increasing sequence of eigenvalues of  problem (\ref{eq 1.1}).

 $(i)$ If $\mu_0\leq \mu<0$ and $N\geq 3$, then  we have that  for $k\in\N$,
 $$\sum^k_{i=1}\lambda_{\mu,i}  \geq  \max\left\{ \Big(1-\frac{\mu}{\mu_0}\Big)   C_N    |\Omega| ^{-\frac{2}{N}} k^{1+\frac{2}{N}}+\frac{\mu}{\mu_0}\lambda_{\mu_0,1}  k,\ b_k e^{-1}\sigma_\mu k^{1+\frac{2}{N}}\right\},$$
 where $C_N$ is the Li-Yau's constant defined in (\ref{LY const}),
 $ b_k >\frac12,\ \lim_{k\to+\infty}b_k= \frac{N}{N+2}$  and
 \begin{equation}\label{1.1-1}
 \sigma_\mu=\max\Big\{ S_N(N-2)^{-\frac{2(N-1)}{N}} \|\mathbb{X}_1^{\frac{1- N}2}\|_{L^2(\Omega)}^{-\frac4N},\,  \big(1-\frac{\mu}{\mu_0}\big)  \frac{N(N-2)}{4}\omega_{_{N-1}}^{\frac2N}|\Omega|^{-\frac2N}\Big\};
 \end{equation}

$(ii)$ If  $\mu>0$ and $N\geq 2$, then  we have that for $k\in\N$,
 $$\sum^k_{i=1}\lambda_{\mu,i}  \geq   C_N   |\Omega|^{-\frac{2}{N}} k^{1+\frac{2}{N}}+ D_0^{-2}\mu\, k.$$

 \end{theorem}

We remark that   the regularity of $\Omega$ can be released to `bounded'  for the lower bounds,   and for $\mu_0\leq \mu<0$, the lower bounds in part $(i)$ consists of the Berezin-Li-Yau's type bounds
and Karachalios's type bounds.  In fact,  the Li-Yau's method can't be applied for $\mu=\mu_0$ because of
the factor $\big(1-\frac{\mu}{\mu_0}\big)$, which vanishes as $\mu\to \mu_0^+$. This also occurs for the Karachalios's estimate (\ref{test 1-0}).
 To overcome this decay,  we develop the Karacholios's method and our lower bound appears
$e^{-1}\sigma_\mu k^{1+\frac{2}{N}}$, where $\sigma_\mu$ has the nondecreasing monotonicity for $\mu\in[\mu_0,0)$; Particularly
$$\sigma_\mu\geq \sigma_{\mu_0}=S_N(N-2)^{-\frac{2(N-1)}{N}} \|\mathbb{X}_1^{\frac{1- N}2}\|_{L^2(\Omega)}^{-\frac4N}>0.$$

 From the monotonicity of $\{\lambda_{\mu,i}\}_{i\in\N}$, we have the following corollary:

\begin{corollary}\label{cr 1.2}
Let  $\{\lambda_{\mu,i}\}_{i\in\N}$ be   the increasing sequence of eigenvalues of  problem (\ref{eq 1.1}).

 $(i)$ If $\mu_0\leq \mu<0$  and $N\geq 3$, then  for  $k\in\N$,
 $$ \lambda_{\mu,k} \geq  \max\left\{ \Big(1-\frac{\mu}{\mu_0}\Big) C_N     |\Omega|^{-\frac{2}{N}} k^{\frac{2}{N}}+\frac{\mu}{\mu_0}\lambda_{\mu_0,1},\ b_k e^{-1}\sigma_\mu k^{\frac{2}{N}}\right\}.$$

  $(ii)$ If $\mu>0$ and $N\geq 2$,  then for $k\in\N$,
 $$ \lambda_{\mu,k}  \geq   C_N    |\Omega|^{-\frac{2}{N}} k^{\frac{2}{N}}+ D_0^{-2}\mu.$$

 \end{corollary}

In order to obtain  upper bounds, we extend Cheng-Yang's type inequality for $\mu\geq \mu_0$ and our upper bounds  state as follows:

 \begin{theorem}\label{teo u 1.1}  Let $\{\lambda_{\mu,i}\}_{i\in\N}$  be the increasing sequence of eigenvalues of  problem (\ref{eq 1.1}).

 $(i)$  If  $\mu\geq 0$ and $N\geq 2$, then
\begin{equation}\label{chengyang-1}
\lambda_{\mu,k}\leq \Big(1+\frac4N\Big)k^{\frac2N}\lambda_{\mu,1}.
\end{equation}

$(ii)$ If  $\mu\in[\mu_0,0)$ and $N\geq 3$,  then
\begin{equation}\label{chengyang-2}
\lambda_{\mu,k}\leq \Big(1+\frac4N\Big)k^{\frac2N}\lambda_{0,1}+\mu D_0^{-2}.
\end{equation}
 \end{theorem}

 It is worth noting that   the crucial inequality for obtaining the Cheng-Yang's type  inequality (\ref{chengyang-1}) is the following
 \begin{equation}\label{yang0}
\sum^k_{i=1}(\lambda_{\mu,k+1}-\lambda_{\mu,i})^2\leq  2a_\mu \sum^k_{i=1} (\lambda_{\mu,k+1}-\lambda_{\mu,i})\lambda_{\mu,i},
\end{equation}
which holds  with $a_\mu =\frac{2 }{N}$ if $\mu>0$.
However, (\ref{yang0}) has the constant $a_\mu=\frac{2 }{N}\frac{-\mu_0}{\mu-\mu_0}$ if $\mu\in(\mu_0,0)$ for $N\geq 3$, then it will give an equality like
$\lambda_{\mu,k}\leq \Big(1+a_\mu\Big)k^{ a_\mu }\lambda_{\mu,1}$, which is too rough since $a_\mu>\frac2N$
and $\displaystyle \lim_{\mu\to \mu_0^+} a_\mu = +\infty$. So we adopt an approach
of comparing with the Laplacian directly.

Combining the estimate of the first eigenvalue $\lambda_{\mu,1}$, we have the following corollary:
 \begin{corollary}\label{cr 1.1}
For $N\geq 3$ and $\mu>0$, we have that
\begin{equation}\label{chy-1}
\lambda_{\mu,k}\leq \big(1+\frac4N\big)\Big(\lambda_{0,1}+\mu \|\phi_{0,1}\|_{L^\infty}^2\int_{\Omega}|x|^{-2}dx\Big)  k^{\frac2N}.
\end{equation}
For $N\geq 2$ and $\mu>0$,
\begin{equation}\label{chy-2}
\lambda_{\mu,k}\leq \big(1+\frac4N\big)\bigg(\lambda_{0,1} +c_1^{-2}\tau_+(\mu)\|\phi_{0,1}\|^2_{C^1}  \frac{\displaystyle\int_{\Omega}  |x|^{2\tau_+(\mu)-2} dx }{\displaystyle\int_{\Omega}\rho^2(x)  |x|^{2\tau_+(\mu)}  dx }
  \bigg)
k^{\frac2N},
\end{equation}
where $c_0>0$, $\tau_+(\mu)=  -\frac{ N-2}2+\sqrt{\mu-\mu_0}$, $(\lambda_{0,1},\phi_{0,1})$ is the first eigenvalue and associated eigenfunctions of (\ref{eq 1.1}) for $\mu=0$ and $\rho(x)={\rm dist}(x,\partial\Omega)$.

 \end{corollary}

 It is remarkable that the upper bound (\ref{chy-1}) is much simple, but it can't be used in the case that $N=2$ and the upper bound (\ref{chy-2}) is available for $N\geq 2$ with the help of the value $\tau_+(\mu)>0$ for $\mu>0$.

We summarize the
asymptotic of the sum of the first $k$-eigenvalues from Corollary \ref{cr 1.2} and Theorem \ref{teo u 1.1} in  the following table:
\begin{center}
\renewcommand{\arraystretch}{1.6}
    \begin{tabular}{ | l | l | l |}
    \hline
    $\mu$ & $(0, +\infty)$ &  $[\mu_0,0)$  \\ \hline
    LB of  $\lambda_{\mu,k}$ & $C_N |\Omega|^{-\frac{2}{N}}k^{\frac{2}{N}}+ \mu D_0^{-2}$  &   $\max\left\{ \Big(1-\frac{\mu}{\mu_0}\Big) C_N  |\Omega|^{-\frac{2}{N}} k^{\frac{2}{N}}+\frac{\mu}{\mu_0}\lambda_{\mu,1} \Big),\, b_k e^{-1}\sigma_\mu k^{\frac{2}{N}}\right\} $
\\ \hline
UB of  $\lambda_{\mu,k}$ & $\Big(1+\frac4N\Big)k^{\frac2N}\lambda_{\mu,1}$ & $ \Big(1+\frac4N\Big)k^{\frac2N}\lambda_{0,1}+\mu D_0^{-2}$  \\
\hline
\end{tabular}\\[1.5mm]
where  $LB$ and $UB$ stand for   Lower bound and Upper bound respectively. \qquad \qquad \qquad
 \qquad  \renewcommand{\arraystretch}{1}
\end{center}

Finally, we provide the Weyl's limit of eigenvalues for Hardy-Leray operators.

\begin{theorem}\label{teo limit}
 Assume that $N\geq 2$, $ \mu\geq \mu_0$ and  $\{\lambda_{\mu,i}\}_{i\in\N}$  is the increasing sequence of eigenvalues of  problem (\ref{eq 1.1}).   Then there holds
\begin{equation}\label{lim 1}
\lim_{k\to+\infty}\lambda_{\mu,k}k^{-\frac{2}{N}}= c_N     |\Omega|   ^{-\frac{2}{N}},
\end{equation}
where $c_N$ is the Weyl's constant given in (\ref{W const}).

 \end{theorem}

We notice that  the limit of $\{\lambda_{\mu,k}k^{-\frac{2}{N}}\}_{k\in\N}$ is independent of $\mu$. Actually, this limit coincides the one for Laplacian, see (\ref{W const}). This answers our question that the   inverse square potential term is a second role for the asymptotic eigenvalues. \smallskip

The rest of this paper is organized as follows. In Section 2, we  build the Berezin-Li-Yau's type lower bounds and
the Karachalios's type lower bounds in Theorem \ref{teo 1.2}.  Section 3 is devoted to the Yang's inequality   and proof of Theorem \ref{teo u 1.1}.   Finally, we prove the Weyl's limit of eigenvalues in Theorem \ref{teo limit} in Section 4.


 \section{Lower bounds}\label{ }

\subsection{Berezin-Li-Yau's type lower bounds}

\begin{proposition}\label{pr 3.1}
Assume that $N\geq3$, $\mu_0<\mu<0$ and   $\{\lambda_{\mu,i} \}_{i\in\N}$  is the increasing sequence of eigenvalues of  problem (\ref{eq 1.1}).
  Then we have

$(i)$ for $\mu_0<\mu<0$ and $k\in\N$ there holds
 $$\sum^k_{i=1}\lambda_{\mu,i}  \geq   \Big(1-\frac{\mu}{\mu_0}\Big)  C_N      |\Omega|^{-\frac{2}{N}} k^{1+\frac{2}{N}}+\frac{\mu}{\mu_0}\lambda_{\mu_0,1}  k;$$

 $(ii)$ for $\mu>0$ and $k\in\N$ there holds
 $$\sum^k_{i=1}\lambda_{\mu,i}  \geq    C_N    |\Omega|^{-\frac{2}{N}} k^{1+\frac{2}{N}}+ D_0^{-2}\mu\, k ,$$
 where we recall  $D_0=\max_{x\in\partial\Omega}|x|$.
 \end{proposition}

For $N\geq 2$, $\mu\geq \mu_0$, we denote
$H_{\mu}(\Omega)$ as the completion of $C_0^\infty(\Omega)$ with the norm
$$\|u\|_{\mu}=\sqrt{\int_\Omega |\nabla u |^2  dx +\mu\int_\Omega \frac{u^2}{|x|^2}dx} $$
and it is a Hilbert space with the inner product
$$\langle u,v\rangle_{\mu}=\int_\Omega  \nabla u \cdot\nabla v   dx +\mu\int_\Omega \frac{u v }{|x|^2}dx. $$
We remark that
$H_{\mu}(\Omega)=H^1_0(\Omega)$  for $\mu>\mu_0$ and if $\mu_0\not=0$, $H_{\mu_0}(\Omega) \varsupsetneqq H^1_0(\Omega).$
 The following lemma is  crucial to get the Berezin-Li-Yau's lower bound whch is appeared in \cite{LY}.

\begin{lemma}\label{lm 3.1}\cite[Lemma 1]{LY} If $f$ is a real-valued function defined on $\R^N$ with $0\leq f\leq M_1$ and
$$\int_{\R^N} f(z)|z|^2 dz\leq M_2,$$
then we have
$$\int_{\R^N} f(z) dz\leq ( M_1|B_1|)^{\frac{2}{N+2}}\Big(M_2\frac{N+2}{N}\Big)^{\frac{N}{N+2}}.$$
\end{lemma}

\noindent{\bf Proof of Proposition \ref{pr 3.1}. }
Let
$(\lambda_{\mu,k}, \phi_k)$ be the eigenvalue and eigenfunction pair of (\ref{eq 1.1})
 such that for positive integer $k$
 $$\|\phi_k\|_{L^2(\Omega)}=1.$$
Then
$$\lambda_{\mu,k} =\min \Big\{\int_{\Omega} \Big(|\nabla u|^2+\frac{ \mu}{|x|^2} u^2\Big)dx \::\:  u\in \cH_k(\Omega)\: \text{with} \: \norm{u}_{L^2(\Omega)}=1\Big\},$$ where
$$
\cH_1(\Omega)= H_\mu(\Omega)\quad \text{and}\quad  \cH_k(\Omega)=\{u\in H_\mu(\Omega)\::\: \text{$\int_{\Omega} u \phi_j \,dx =0$ for $j=1,\dots, k-1$}\}\quad \text{for $k>1$.}
$$
Moreover, $\{\phi_k \in \cH_\mu (\Omega):\: k \in \N\}$ is an orthonormal basis of $L^2(\Omega)$.

Denote
$$\Phi_k(x,y):=\sum^k_{j=1} \phi_j(x)\phi_j(y)$$
and its  Fourier transform is then
$$\hat{\Phi}_k(z,y)=(2\pi)^{-\frac N2}\int_{\R^N} \Phi_k(x,y) e^{-{\rm i}\,x\cdot z} dx. $$
Note that
\begin{align*}
 \int_{ \R^N}  \int_{\R^N} |\hat{\Phi}_k(z,y)|^2dzdy  =  \int_{\Omega} \int_{\R^N}|\Phi_k(x,y)|^2dxdy
 =  \sum^k_{j=1} \int_{\Omega} \phi^2_j(x) dx=k
 \end{align*}
 and
  \begin{align*}
   \int_{\R^N} |\hat{\Phi}_k(z,y)|^2dy  &=(2\pi)^{-  N }  \int_{\Omega}  \big|\int_{\Omega}\Phi_k(x,y)e^{-{\rm i}\,x\cdot z}dx \big|^2  dy
   =(2\pi)^{-  N }\big|\int_{\Omega}\sum^k_{j=1}\phi_j(x)e^{-{\rm i}\, x\cdot z}dx \big|^2,
 \end{align*}
 which implies by Bessel's inequality (see \cite[(1.2)]{MF}) that
\begin{align*}
   \int_{\R^N} |\hat{\Phi}_k(z,y)|^2dy  &\leq (2\pi)^{-  N }    \int_{\Omega}|e^{{\rm -i}\, x\cdot z}|^2 dx =(2\pi)^{-  N } |\Omega|.
 \end{align*}
Meanwhile, the Hardy inequality (\ref{Hardy inequalities})  implies that  for $\mu_0<\mu<0$,
 \begin{align*}
 \int_\Omega |\nabla u |^2  dx+\mu\int_\Omega \frac{u^2 }{|x|^2}dx&=\frac{\mu}{\mu_0}\Big(\int_\Omega |\nabla u |^2  dx +\mu_0\int_\Omega \frac{u^2 }{|x|^2}dx\Big)
 \\&\quad +(1-\frac{\mu}{\mu_0})\int_\Omega |\nabla u |^2  dx
 \\&\geq (1-\frac{\mu}{\mu_0})\int_\Omega |\nabla u |^2  dx+\frac{\mu}{\mu_0}\lambda_{\mu_0,1} \int_\Omega  u^2 dx,
 \end{align*}
thus, we have  that
 \begin{align*}
   \int_{\R^N}\int_{\R^N} &|\hat{\Phi}_k(z,y)|^2 |z|^2dydz   =\int_{\R^N}\int_\Omega |\nabla \Phi_k(x,y) |^2 dy dx
    =\int_{\Omega} \Big|\sum^k_{j=1}\nabla\phi_j(x)\Big|^2 dx
   \\&\leq (1-\frac{\mu}{\mu_0})^{-1}\Big(\int_\Omega \Big|\nabla\sum^k_{j=1}\phi_j(x) \Big|^2  dx+\mu\int_\Omega \frac{\sum^k_{j=1} \phi_j^2(x)}{|x|^2}dx-\frac{\mu}{\mu_0}\lambda_{\mu_0,1}\int_\Omega  \sum^k_{j=1} \phi_j^2(x) dx\Big)
   \\&=(1-\frac{\mu}{\mu_0})^{-1}\Big(\sum^k_{j=1}\lambda_{\mu,j} -\frac{\mu}{\mu_0}\lambda_{\mu_0,1} k\Big).
 \end{align*}

 For   $ \mu>0$, we have that
 \begin{align*}
 \int_\Omega |\nabla u |^2  dx+\mu\int_\Omega \frac{u^2 }{|x|^2}dx
 \geq  \int_\Omega |\nabla u |^2  dx+ \mu D_0^{-2}  \int_\Omega  u^2 dx,
 \end{align*}
and then
\begin{align*}
   \int_{\R^N}\int_{\R^N}  |\hat{\Phi}_k(z,y)|^2 |z|^2dydz   &=\int_{\R^N}\int_\Omega |\nabla \Phi_k(x,y) |^2 dy dx
   \\&\leq \int_\Omega \Big| \sum^k_{j=1}\nabla\phi_j \Big|^2  dx+\mu\int_\Omega \frac{\sum^k_{j=1}\phi_j^2 }{|x|^2}dx-D_0^{-2}\mu\int_\Omega  \Big(\sum^k_{j=1}\phi_j\Big)^2 dx
   \\&=  \sum^k_{j=1}\lambda_{\mu,j} -\mu D_0^{-2}    k.
 \end{align*}

Now for {\it Case: $N\geq 3$, $\mu_0<\mu<0$.}  We apply Lemma \ref{lm 3.1} to the function $f(z)=\displaystyle \int_{\Omega} |\hat{\Phi}_k(z,y)|^2dy$
 with
 $$M_1=(2\pi)^{-  N} |\Omega| \quad{\rm and}\quad M_2=(1-\frac{\mu}{\mu_0})^{-1} \Big(\sum^k_{j=1}\lambda_{\mu,j} -\frac{\mu}{\mu_0}\lambda_{\mu_0,1}  k\Big),$$
  then we conclude that
 \begin{align*}
k = \int_{\R^N} f(z) dz  &\leq ( M_1|B_1|)^{\frac{2}{N+2}} \Big(M_2\frac{N+2}{N}\Big)^{\frac{N}{N+2}}
  \\&\leq \Big( (2\pi)^{- N} |\Omega|\,  |B_1|\Big)^{\frac{2}{N+2}}\Big[ (1-\frac{\mu}{\mu_0})^{-1} \Big(\sum^k_{j=1}\lambda_{\mu,j} -\frac{\mu}{\mu_0}\lambda_{\mu_0,1} k\Big)\Big]^{\frac{N}{N+2}}\Big(\frac{N+2}{N}\Big)^{\frac{N}{N+2}}
 \end{align*}
and
$$\sum^k_{j=1}\lambda_{\mu,j} \geq \Big(1-\frac{\mu}{\mu_0}\Big) C_N  |\Omega| ^{-\frac{2}{N}} k^{1+\frac{2}{N}}+\frac{\mu}{\mu_0}\lambda_{\mu_0,1}  k .$$

For {\it Case: $\mu>0$.}  We again apply Lemma \ref{lm 3.1} to the same function $f(z)$
and $M_1$, but $$M_2=   \sum^k_{j=1}\lambda_{\mu,j}-\mu D_0^{-2}    k ,$$
  then we conclude that
 \begin{align*}
k&= \int_{\R^N} f(z) dz  \leq ( M_1|B_1|)^{\frac{2}{N+2}}M_2^{\frac{N}{N+2}}\Big(\frac{N+2}{N}\Big)^{\frac{N}{N+2}}
  \\&\leq \Big( (2\pi)^{-  N } |\Omega|\,  |B_1|\Big)^{\frac{2}{N+2}} \Big(\sum^k_{j=1}\lambda_{\mu,j} -\mu D_0^{-2} k\Big) ^{\frac{N}{N+2}}\Big(\frac{N+2}{N}\Big)^{\frac{N}{N+2}}
 \end{align*}
and
$$\sum^k_{j=1}\lambda_{\mu,j} \geq  C_N   |\Omega| ^{-\frac{2}{N}} k^{1+\frac{2}{N}}+ \mu D_0^{-2} k .$$
We complete the proof. \hfill$\Box$\medskip

 \begin{corollary}\label{cr 3.1}
Let  $\{\lambda_{\mu,i} \}_{i\in\N}$ be the increasing sequence of eigenvalues of  problem (\ref{eq 1.1}).

  $(i)$ For $N\geq 3$, $\mu_0<\mu<0$ and $k\in\N$, we have that
 $$ \lambda_{\mu,k} \geq   \Big(1-\frac{\mu}{\mu_0}\Big) \Big(C_N |\Omega|^{-\frac{2}{N}} k^{\frac{2}{N}}+\frac{\mu}{\mu_0}\lambda_{\mu_0,1} \Big). $$

 $(ii)$ For $N\geq 2$, $\mu>0$ and $k\in\N$,  we have that
 $$ \lambda_{\mu,k} \geq   C_N |\Omega| ^{-\frac{2}{N}} k^{\frac{2}{N}}+\mu D_0^{-2} .$$
 \end{corollary}
\noindent{\bf Proof. }   From the increasing monotonicity of $k\mapsto \lambda_{\mu,k}$,
we have that
  $$\lambda_{\mu,k}\geq \frac1k\sum^k_{j=1}\lambda_{\mu,j},$$
  which implies the lower bounds for $\lambda_{\mu,k}$  by  Proposition \ref{pr 3.1}.
  \hfill$\Box$

\subsection{Karachalios' type lower bound}

Recall  that $S_N=2^{2/N}\pi^{1+1/N}\Gamma(\frac{N+1}{2})$ is the best constant of Sobolev inequality in $\R^N$ and $\Gamma$ is the Gamma function. In this subsection, we would use the following notations
$$\mathbb{X}_1(x)=\big(-\log(\frac{|x|}{D_0})\big)^{-1}\ \ {\rm and}\ \  \mathbb{X}_2(x)\equiv 1,\quad \forall\, x\in\Omega,$$
where we recall $D_0=\max_{x\in\partial\Omega}|x|$.  Note that
\begin{align*}
\|\mathbb{X}_1^{-\frac{ N-1}{2}}\|_{L^2(\Omega)}^2  &\leq  \int_{B_{D_0}}  \mathbb{X}_1^{1- N}(y) dy
  \\&\leq  \omega_{_{N-1}} \int_0^{D_0} r^{N-1} (-\log \frac{r}{D_0})^{N-1} dr < +\infty
   \end{align*}
and $\|\mathbb{X}_2^{-\frac{ N-1}{2}}\|_{L^2(\Omega)}^2  =|\Omega|.$

\begin{proposition}\label{pr 3.2}
Assume that $N\geq 3$,  $\mu_0 \leq \mu<0$ and  $\{\lambda_{\mu,k} \}_{k\in\N}$  is the increasing sequence of eigenvalues of  problem (\ref{eq 1.1}).
 Then for $k\in\N$ we have that
 $$ \lambda_{\mu,k}  \geq  e^{-1}\sigma_\mu k^{\frac{2}{N}},$$
 where  $\sigma_\mu$ is defined in (\ref{1.1-1}).

 \end{proposition}
\noindent{\bf Proof. }  When  $\mu=\mu_0$, the proof is refereed to \cite{K2} and we focus on the case $\mu_0<\mu<0$.

 Let  $\cL_\mu$ be the Hardy-Leray operator with its domain defined as
$D(\cL_\mu) =C_0^\infty(\Omega)$ and its Friedrich's extension, still denoted by $\cL_\mu$, with its domain defined as
$$D(\cL_\mu):=\big\{u\in H_\mu(\Omega): \cL_\mu u\in L^2(\Omega)\big\},$$
which is a nonnegative self-adjoint operator on $L^2(\Omega)$ and the operator gives rise to the semigroup  of operators $e^{-\cL_\mu t}$ for every $t>0$, possessing an integral kernel $K(x,y,t)>0$ for all $(x,y,t)\in \Omega\times\Omega\times(0,+\infty)$. Then $\cL_\mu$ has compact resolvent and $K$ can be represented as
$$K(x,y,t)=\sum^\infty_{i=1}e^{-\lambda_{\mu,i}  t} \phi_i(x)\phi_i(y),$$
which solves the problem
\begin{equation}\label{eq 3.1}
\left\{ \arraycolsep=1pt
\begin{array}{lll}
\partial_t K+ \cL_\mu K=0\quad \  &{\rm in}\ \,   \Omega\times\Omega\times(0,+\infty),\\[2mm]
 \phantom{ \partial_t \,  }
  K(x,y,t)>0\quad \ &{\rm{in}}\ \,   \Omega\times\Omega\times(0,+\infty),\\[2mm]
 \phantom{ \partial_t  \, }
  K(x,y,t)=0\quad \ &{\rm{on}}\ \,   \partial\Omega\times\partial\Omega\times(0,+\infty).
\end{array}
\right.
\end{equation}
Since $\{\phi_i\}_{i\geq 1}$ is an orthonormal basis of $L^2(\Omega)$, it follows that
$$h(t):=\sum^\infty_{i=1} e^{-2\lambda_{\mu,i} t}=\int_\Omega\int_\Omega K^2(x,y,t) dxdy. $$

By the improved Hardy-Sobolev inequalities, \cite[Theorem A]{FT}
for $\mu\geq \mu_0$, (the following inequality is sharp for $\mu= \mu_0$)
$$\int_\Omega |\nabla u|^2 dx +\mu\int_\Omega \frac{u^2}{|x|^2} dx\geq \sigma_1 \Big(\int_{\Omega} |u|^{2^*} \mathbb{X}_1(x)^{-\frac{2(N-1)}{N-2}} dx\Big)^{\frac{2}{2^*}},\quad \forall\, u\in C_0^\infty(\Omega)$$
and for $\mu>\mu_0$
$$\int_\Omega |\nabla u|^2 dx +\mu  \int_\Omega \frac{u^2}{|x|^2} dx\geq \sigma_2 \Big(\int_{\Omega} |u|^{2^*}  dx\Big)^{\frac{2}{2^*}},\quad \forall\, u\in C_0^\infty(\Omega),$$
where $2^*=\frac{2N}{N-2}$,
$$\sigma_1=S_N(N-2)^{-\frac{2(N-1)}{N}} \quad
{\rm and} \quad \sigma_2=\frac{N(N-2)}{4}(1-\frac{\mu}{\mu_0}) \omega_{_{N-1}}.$$
 By H\"older inequality, we get that for $j=1,2$,
 \begin{align}
h(t)&\leq   \int_\Omega\Big(\big(\int_\Omega \mathbb{X}_j(y)^{\frac{2(N-1)}{N-2}} K^{2^*}(x,y,t) dy\big)^{\frac1{2^*-1}}\big(\int_\Omega \mathbb{X}_j(y)^{-\frac{2(N-1)}{(N-2)(2^*-2)}} K(x,y,t) dy\big)^{\frac{2^*-2}{2^*-1}}\Big)dx\nonumber
  \\&\leq  \Big[\int_\Omega\Big(\int_\Omega \mathbb{X}_j(y)^{\frac{2(N-1)}{N-2}} K^{2^*}(x,y,t) dy\Big)^{\frac2{2^*}} dx\Big]^{\frac{2^*}{2(2^*-1)}}\Big(\int_\Omega Q_j^2(x,t) dx\Big)^{\frac{2^*-2}{2(2^*-1)}}, \label{ineq 1-1}
 \end{align}
 where $\frac{2(N-1)}{(N-2)(2^*-2)}=\frac{ N-1}{2}$ and
$$Q_j(x,t)=\int_\Omega \mathbb{X}_j(y)^{-\frac{ N-1}{2}} K(x,y,t) dy.  $$

We observe that $Q_j(x,t)$ is the solution of the Cauchy-Dirichlet problem:
\begin{equation}\label{eq 3.2}
\left\{ \arraycolsep=1pt
\begin{array}{lll}
\partial_t Q_j+ \cL_\mu Q_j=0\quad \  &{\rm in}\ \   \Omega\times(0,+\infty),\\[2mm]
 \phantom{ \partial_t \, \partial_t \ }
  Q_j(x, 0)=\mathbb{X}_j^{\frac{1- N}2}\quad \ &{\rm{in}}\ \  \Omega,\\[2mm]
 \phantom{ \partial_t \, \partial_t\ \, }
  Q_j(x, t)=0\quad \ &{\rm{on}}\ \  \partial\Omega\times(0,+\infty).
\end{array}
\right.
\end{equation}
Multiplying the above equation (\ref{eq 3.2}) by $Q_j$, we get the energy equation
$$\frac12\frac{d}{dt}\|Q_j(t)\|_{L^2(\Omega)}^2+\|Q_j(t)\|^2_{\mu}=0$$
and then
$$\|Q_j(t)\|_{L^2(\Omega)}^2\leq \|Q_j(0)\|_{L^2(\Omega)}^2=\| \mathbb{X}_j^{\frac{1- N}2}\|_{L^2(\Omega)}^2.$$

With the help of above inequalities,   letting
$$C_j:=  \|\mathbb{X}_j^{\frac{1- N}2}\|_{L^2(\Omega)}^{\frac{2(2^*-1)}{2^*}}\quad{\rm for}\ \ j=1,2,$$
 for $\mu>\mu_0$,  from (\ref{ineq 1-1})  we obtain that
\begin{align*} h^{\frac{2(2^*-1)}{2^*}}(t)&\leq  C_j \int_\Omega\Big(\int_\Omega \mathbb{X}_j(y)^{\frac{2(N-1)}{N-2}} K^{2^*}(x,y,t) dy\Big)^{\frac2{2^*}} dx
\\&\leq  \frac{C_j}{\sigma_j}
\int_\Omega\int_\Omega\Big(  |\nabla_y K(x,y,t)|^2 dx +\mu \frac{K^2(x,y,t)}{|y|^2} \Big)dxdy
\\&\leq \frac{C_j}{\sigma_j}\|Q_j(t)\|^2_{\mu} =-\frac{C_j}{\sigma_j} \frac12\frac{d}{dt}\|Q_j(t)\|_{L^2(\Omega)}^2
\\&\leq
-\frac{C_j}{\sigma_j} \frac12\frac{dh(t)}{dt},
\end{align*}
which, using $h(0)=+\infty$, implies that
$$\sum^\infty_{i=1} e^{-2\lambda_{\mu,i} t}=h(t)\leq \Big(\frac{2^*}{2\sigma_j (2^*-2) }\Big)^{\frac{2^*}{2^*-2}}\|\mathbb{X}_j^{\frac{1- N}2}\|_{L^2(\Omega)}^2 t^{-\frac{2^*}{2^*-2}}\quad{\rm for}\ \ j=1,2.$$
Now we choose
$$t=\frac{2^*}{2(2^*-2)}\frac1{\lambda_{\mu,k} }$$
and
$$ke^{-\frac{2^*}{2^*-2}} \leq \sum^\infty_{i=1} e^{-\frac{2^*\lambda_{\mu,i} }{\lambda_{\mu,k} (2^*-2)}}\leq \sigma_j^{-\frac{2^*}{2^*-2}} \lambda_{\mu,k}  ^{\frac{2^*}{2^*-2}} \|\mathbb{X}_j^{\frac{1- N}2}\|_{L^2(\Omega)}^2\quad{\rm for}\ \ j=1,2. $$
Note that $\frac{2^*}{2^*-2}=\frac N2$ and we conclude that
$$\lambda_{\mu,k} \geq e^{-1}\sigma_j\|\mathbb{X}_j^{\frac{1- N}2}\|_{L^2(\Omega)}^{-\frac4N} k^{\frac2N}\quad{\rm for }\ j=1,2.$$
As a consequence, we obtain that
$$\lambda_{\mu,k} \geq e^{-1}\sigma_\mu k^{\frac2N},$$
where
$$\sigma_\mu= \max\Big\{  \sigma_1\|\mathbb{X}_1^{\frac{1- N}2}\|_{L^2(\Omega)}^{-\frac4N},\,  \sigma_2|\Omega|^{-\frac2N}\Big\}.$$
We complete the proof. \hfill$\Box$\medskip

\noindent{\bf Proof of Theorem \ref{teo 1.2}. } When $N\geq 3$ and $\mu_0\leq \mu<0$, the lower bounds
$$\sum^k_{i=1}\lambda_{\mu,i} \geq    \Big(1-\frac{\mu}{\mu_0}\Big) \Big( C_N |\Omega| ^{-\frac{2}{N}} k^{1+\frac{2}{N}}+\frac{\mu}{\mu_0}\lambda_{\mu_0,1} k\Big)$$
follows from Proposition \ref{pr 3.1} part $(i)$. Note that for  $\mu=\mu_0$, the above bounds is true obviously since $1-\frac{\mu}{\mu_0}=0$.
Moreover, from Proposition \ref{pr 3.2}, we have that
$$\lambda_{\mu,k} \geq e^{-1}\sigma_\mu k^{\frac2N},$$
which  implies that
$$\sum^k_{i=1}\lambda_{\mu,i}(\Omega)\geq e^{-1}\sigma_\mu k^{\frac 2N}\sum^k_{i=1}\big(\frac{i}{k}\big)^\frac2{N} = b_k e^{-1}\sigma_\mu k^{\frac2N+1},$$
where
\begin{align*}
b_k=\frac1k\sum^k_{i=1}\big(\frac{i}{k}\big)^\frac2{N}\to \frac{N}{N+2} \quad{\rm as}\ \, k\to+\infty
\end{align*}
and for any $k\geq 2$
\begin{align*}
b_k=\frac1k \sum^k_{i=1}\big(\frac{i}{k}\big)^\frac2{N} \geq \frac1k \sum^k_{i=1}\frac{i}{k} =   \frac{k+1}{2k} >  \frac{1}{2}.
\end{align*}
When   $\mu>0$, it comes from Proposition \ref{pr 3.1} part $(ii)$ directly and we complete the proof.\hfill$\Box$\medskip

\noindent{\bf Proof of Corollary \ref{cr 1.2}. }
It follows form Corollary \ref{cr 3.1} and
Proposition \ref{pr 3.2} directly. \hfill$\Box$

\section{Upper bounds }

Now we establish the upper bounds, more precisely we will use the following Cheng-Yang's type inequality and extend it for $\mu \geq \mu_0$.
\begin{proposition}\label{pr yang}\cite[Proposition 3.1]{CZ}
Let $0<\nu_1\leq \nu_2\leq \cdots$  be a sequence of numbers satisfying the
following inequality
\begin{equation}\label{yang1}
 \sum^k_{i=1}(\nu_{k+1}-\nu_i)^2\leq 2\varrho_0\sum^k_{i=1}(\nu_{k+1}-\nu_i)\nu_i,
 \end{equation}
where $\varrho_0$  is a positive constant. Then we have
$$\nu_{k+1}  \leq (1 + 2 \varrho_0)  k^{ \varrho_0} \nu_1.$$
\end{proposition}

From Proposition (\ref{pr yang}), the key point is to obtain inequality (\ref{yang1}) with
explicit constant $\varrho_0$. To this end, we have the following inequalities.

\begin{proposition}\label{pr 2.1}
Assume that $\mu>\mu_0$  and let  $\{\lambda_{\mu,i}\}_{i\in\N}$  be the increasing sequence of eigenvalues of  problem (\ref{eq 1.1}).

 Then $(i)$ for $N\geq 2$ and $\mu>0$, one has
\begin{equation}\label{yang 1-1}
\sum^k_{i=1}(\lambda_{\mu,k+1}-\lambda_{\mu,i})^2\leq \frac{4 }{N}\sum^k_{i=1} (\lambda_{\mu,k+1}-\lambda_{\mu,i})\lambda_{\mu,i};
\end{equation}

$(ii)$ for $N\geq 3$ and $\mu_0<\mu<0$, one has
\begin{equation}\label{yang 1-2}
\sum^k_{i=1}(\lambda_{\mu,k+1}-\lambda_{\mu,i})^2\leq \frac{4 }{N}\frac{-\mu_0}{\mu-\mu_0}\sum^k_{i=1} (\lambda_{\mu,k+1}-\lambda_{\mu,i})\lambda_{\mu,i}.
\end{equation}
 \end{proposition}
\noindent{\bf Proof. } The proof for $\mu=0$  can be found in \cite[Theorem 2.1]{CgY}.
Here the difference is to deal with the Hardy term and we give the proof for reader's convenience.

Let $(\lambda_{\mu,i},\phi_i)$ be the $i$-th eigenvalue and eigenfunctions of (\ref{eq 1.1})
with $\displaystyle\int_{\Omega} \phi_i\phi_jdx=\delta_{ij},$
where $\delta_{ij}=0$ if $i\not=j$ and $\delta_{ii}=1$.
Let $x=(x^1,\cdots,x^N)$  and $g=x^m$ with $m=1,\cdots,N$ define a trial function $\varphi_i$ by
  $$\varphi_i:= g\phi_i-\sum^k_{j=1} a_{ij}\phi_j,\quad a_{ij}:=\int_{\Omega} g\phi_i\phi_j dx=a_{ji}.$$
By the orthonormality of $\phi_i$ and the definition of $a_{ij}$, $\varphi_i$ is perpendicular to $\phi_j$, that is, for $i,j=1,\cdots,k$, there holds $\displaystyle \int_{\Omega} \varphi_i \phi_j dx=0$.
Let
$$b_{ij}= \int_{\Omega} \phi_j \nabla g\cdot \nabla \phi_i dx$$
 and then
 \begin{align*}
  \lambda_{\mu,j}a_{ij} =  \int_{\Omega}  g(\cL_\mu \phi_j) \phi_idx
 &=    \int_{\Omega}\Big(-2\phi_j \nabla g\cdot \nabla \phi_i + g(\cL_\mu \phi_i) \phi_jdx\Big)
  = -2b_{ij}+ \lambda_{\mu,i}a_{ij},
 \end{align*}
 i.e.
\begin{equation}\label{e1}
2b_{ij}=(\lambda_{\mu,i}-\lambda_{\mu,j})a_{ij}.
\end{equation}
Note that
$$\cL_\mu\varphi_i=\lambda_{\mu,i}g\phi_i-2\nabla g\cdot \nabla \phi_i-\sum^k_{j=1} a_{ij}\lambda_{\mu,j}\phi_j.$$
 Hence, we infer that
 \begin{align*}
 \|\varphi_i\|_\mu^2&=\lambda_{\mu,i} \int_\Omega \varphi_i^2dx-2\int_\Omega \varphi_i \nabla g\cdot \nabla \phi_i dx\\
\\ &= \lambda_{\mu,i} \int_\Omega \varphi_i^2dx -2  \int_{\Omega} (g \nabla g)\cdot (\phi_i\nabla \phi_i)dx  + 2\sum^k_{j=1} a_{ij} \int_\Omega \phi_j\nabla g\cdot \nabla \phi_i dx
 \\&=\lambda_{\mu,i} \int_\Omega \varphi_i^2dx+1+\sum^k_{j=1} (\lambda_{\mu,i}-\lambda_{\mu,j})a_{ij}^2,
 \end{align*}
 where we note particularly
 \begin{equation}\label{e2}
-2\int_\Omega \varphi_i \nabla g\cdot \nabla \phi_i dx  = 1+\sum^k_{j=1} (\lambda_{\mu,i}-\lambda_{\mu,j})a_{ij}^2.
\end{equation}

From the Rayleigh-Ritz inequality, we have that
$$
 \lambda_{\mu,k+1} \int_\Omega \varphi_i^2dx\leq \|\varphi_i\|_\mu^2=\lambda_{\mu,i} \int_\Omega \varphi_i^2dx+1+\sum^k_{j=1} (\lambda_{\mu,i}-\lambda_{\mu,j})a_{ij}^2,
$$
that is
  \begin{equation}\label{e3}
  (\lambda_{\mu,k+1}-\lambda_{\mu,i} ) \int_\Omega \varphi_i^2dx\leq 1+\sum^k_{j=1} (\lambda_{\mu,i}-\lambda_{\mu,j})a_{ij}^2.
  \end{equation}
Multiplying (\ref{e2}) by $  (\lambda_{\mu,k+1}-\lambda_{\mu,i} )^2$ and then taking sum on $i$ from 1 through $k$, we obtain that
   \begin{align*}
&-2\sum^k_{i=1} (\lambda_{\mu,k+1}-\lambda_{\mu,i} )^2\int_\Omega \varphi_i \nabla g\cdot \nabla \phi_i dx    \\ = &\sum^k_{i=1} (\lambda_{\mu,k+1}-\lambda_{\mu,i} )^2+\sum^k_{i,j=1} (\lambda_{\mu,i}-\lambda_{\mu,j})(\lambda_{\mu,k+1}-\lambda_{\mu,i} )^2a_{ij}^2
\\=& \sum^k_{i=1} (\lambda_{\mu,k+1}-\lambda_{\mu,i} )^2-4\sum^k_{i,j=1} (\lambda_{\mu,k+1}-\lambda_{\mu,i}) b_{ij}^2
\\:=& \theta
 \end{align*}
 by the symmetry of $a_{ij}$ and the anti-symmetry of $b_{ij}$.
 Multiplying (\ref{e3}) by $  (\lambda_{\mu,k+1}-\lambda_{\mu,i} )^2$ and then taking sum on $i$ from 1 through $k$, we obtain that
   \begin{align*}
   \sum^k_{i=1}(\lambda_{\mu,k+1}-\lambda_{\mu,i} )^3 \int_\Omega \varphi_i^2dx\leq  \sum^k_{i=1} (\lambda_{\mu,k+1}-\lambda_{\mu,i} )^2-4\sum^k_{i,j=1} (\lambda_{\mu,k+1}-\lambda_{\mu,i}) b_{ij}^2=\theta.
 \end{align*}
 Note that for arbitrary constant $d_{ij}$ (to be determined later)
 \begin{align*}
   \theta^2&=\bigg(-2\sum^k_{i=1} (\lambda_{\mu,k+1}-\lambda_{\mu,i} )^2\int_\Omega \varphi_i \nabla g\cdot \nabla \phi_i dx \bigg)^2
   \\&=4\bigg( \sum^k_{i=1} \int_\Omega (\lambda_{\mu,k+1}-\lambda_{\mu,i} )^2 \varphi_i \nabla g\cdot \nabla \phi_idx-(\lambda_{\mu,k+1}-\lambda_{\mu,i} )^{\frac32}\int_\Omega \sum^k_{j=1}d_{ij}\varphi_i \phi_j  dx \bigg)^2
  \\&\leq 4\bigg(\sum^k_{i=1} (\lambda_{\mu,k+1}-\lambda_{\mu,i} )^3\int_\Omega\varphi_i^2\bigg)\bigg(\sum^k_{i=1}\int_\Omega\Big((\lambda_{\mu,k+1}-\lambda_{\mu,i} )^{\frac12}\nabla g\cdot \nabla \phi_i-\sum^k_{j=1}d_{ij} \phi_j\Big)^2 dx \bigg)
  \\&\leq 4\theta\! \int_\Omega\Big( \sum^k_{i=1}(\lambda_{\mu,k+1}-\lambda_{\mu,i} ) |  \partial_m \phi_i|^2 -2\sum^k_{i,j=1}d_{ij} (\lambda_{\mu,k+1}-\lambda_{\mu,i} )^{\frac12}\phi_j\nabla g\cdot \nabla \phi_i +(\sum^k_{i,j=1}d_{ij} \phi_j)^2\Big)dx,
 \end{align*}
 then we have that
 \begin{align*}
   \theta &\leq 4  \sum^k_{i=1}(\lambda_{\mu,k+1}-\lambda_{\mu,i} )\int_\Omega |  \partial_m \phi_i|^2dx+4\Big(-2\sum^k_{i,j=1}d_{ij} (\lambda_{\mu,k+1}-\lambda_{\mu,i} )^{\frac12}b_{ij} + \sum^k_{i,j=1}d_{ij}^2 \Big).
 \end{align*}
 Taking $d_{ij}=(\lambda_{\mu,k+1}-\lambda_{\mu,i} )^{\frac12}b_{ij} $,
we have that
   \begin{align*}
& \theta = \sum^k_{i=1} (\lambda_{\mu,k+1}-\lambda_{\mu,i} )^2-4\sum^k_{i,j=1} (\lambda_{\mu,k+1}-\lambda_{\mu,i}) b_{ij}^2
 \\& \leq 4  \sum^k_{i=1}(\lambda_{\mu,k+1}-\lambda_{\mu,i} ) \int_\Omega |  \partial_m \phi_i|^2dx-4\sum^k_{i,j=1} (\lambda_{\mu,k+1}-\lambda_{\mu,i} ) b_{ij}^2.
 \end{align*}
 Thus,
 $$ \sum^k_{i=1} (\lambda_{\mu,k+1}-\lambda_{\mu,i} )^2\leq 4  \sum^k_{i=1}(\lambda_{\mu,k+1}-\lambda_{\mu,i} ) \int_\Omega |  \partial_m \phi_i|^2dx.$$
 Finally, we take sum on $m$ from 1 through $N$ and obtain that
   \begin{align*}
   N \sum^k_{i=1} (\lambda_{\mu,k+1}-\lambda_{\mu,i} )^2&\leq    4\sum^k_{j=1} (\lambda_{\mu,k+1}-\lambda_{\mu,i} )\int_{\Omega}|\nabla \phi_i|^2dx.
  \end{align*}
Therefore we conclude that for $\mu>0$,
   \begin{align*}
   N \sum^k_{i=1} (\lambda_{\mu,k+1}-\lambda_{\mu,i} )^2 & <    4\sum^k_{i=1} (\lambda_{\mu,k+1}-\lambda_{\mu,i} )\Big(\int_{\Omega}|\nabla \phi_i|^2dx+\mu\int_\Omega\frac{\phi_i^2}{|x|^2}dx\Big)
 \\& =4\sum^k_{i=1} (\lambda_{\mu,k+1}-\lambda_{\mu,i} )\lambda_{\mu,i}.
 \end{align*}
 While for $\mu_0<\mu<0$, since we have that
 $$\int_{\Omega}|\nabla \phi_i|^2dx\leq \frac{-\mu_0}{\mu-\mu_0} \Big(\int_{\Omega}(|\nabla \phi_i|^2+\mu\frac{\phi_i^2}{|x|^2})dx\Big)
 =\frac{-\mu_0}{\mu-\mu_0}\lambda_{\mu,i}$$
and then
   \begin{align*}
   N \sum^k_{i=1} (\lambda_{\mu,k+1}-\lambda_{\mu,i} )^2 &\leq = 4\frac{-\mu_0}{\mu-\mu_0} \sum^k_{i=1} (\lambda_{\mu,k+1}-\lambda_{\mu,i} )\lambda_{\mu,i}.
 \end{align*}
We complete the proof. \hfill$\Box$ \medskip

\noindent{\bf Proof of Theorem \ref{teo u 1.1}. }
{\it Part $(i)$: Case of $\mu>0$. }
Theorem \ref{teo u 1.1} follows Proposition \ref{pr yang} and (\ref{yang 1-1}) with $\varrho_0=\frac{2}{N}$. To be convenient,
we sketch the proof here.
Denote
$$\Lambda_k=\frac1k\sum^k_{i=1} \lambda_{\mu,i},\ \ T_k=\frac1k\sum^k_{i=1} \lambda_{\mu,i}^2,\ \ F_k=(1+\frac2N)\Lambda_k^2-T_k,$$
 then we have
   \begin{equation}\label{y 1.1}
   F_{k+1}\leq C(N,k) (\frac{k+1}{k})^{\frac4N} F_k,
   \end{equation}
  where
  $$0< C(N,k):=1-\frac16(\frac{k}{k+1})^{\frac4N} \frac{(1+\frac2N)(1+\frac4N)}{(k+1)^3}<1,$$
 and then
   \begin{align*}
F_k   \leq   C(N,k-1) (\frac{k}{k-1})^{\frac4N} F_{k-1}&<(\frac{k}{k-1})^{\frac4N} F_{k-1}\leq  \cdots
 \\&<  (\frac{k}{k-1})^{\frac4N}(\frac{k-1}{k-2})^{\frac4N}  \cdots (\frac21)^{\frac4N}F_{1}=\frac2N k^{\frac4N} \lambda_{\mu,1}^2.
 \end{align*}
 On the other hand, (\ref{y 1.1}) implies that
   \begin{align*}
   \Big(\lambda_{\mu,k+1} -\frac{2+N}N\Lambda_k\Big)^2\leq  \frac{4+N}NF_k-\frac{4+N}N\frac{2+N}N\Lambda_k^2
 \end{align*}
  and then
  \begin{align*}
 \frac{2}{4+N}\lambda_{\mu,k+1}^2+\frac{4+N}{2+N}\Big(\lambda_{\mu,k+1}-\frac{4+N}N\Lambda_k \Big)^2\leq \frac{4+N}N F_k,
 \end{align*}
  that means
  $$\lambda_{\mu,k+1}^2\leq \frac{N}{N+2}(\frac{N+4}N)^2 F_k\leq (1+\frac4N)^2k^{\frac4N} \lambda_{\mu,1}^2.$$
  Thus, we have that
\begin{equation}\label{e4}
\lambda_{\mu,k}\leq \big(1+\frac4N\big)k^{\frac2N}\lambda_{\mu,1}.
\end{equation}

  \medskip

 {\it Part $(ii)$: Case of $\mu_0\leq \mu<0$. } In this case, the constant $\frac{4 }{N}\frac{-\mu_0}{\mu-\mu_0}>\frac{4 }{N}$, which will produce a higher order for parameter $k$ from Proposition \ref{pr yang}. So we shall  derive the upper bounds  by comparing with the eigenvalues of  Laplacian Dirichlet problem.

  Let
  $$\tilde \cL_\mu =\cL_\mu -\mu D_0^{-2}$$
  and    $\lambda_{\tilde \cL_\mu,k}$ be the $k$-th Dirichlet eigenvalue of
  $\tilde \cL_\mu$. Clearly we have
 $$\lambda_{\tilde \cL_\mu,k}=\lambda_{\mu,k}-\mu D_0^{-2}.$$
  Then for $\mu\in[\mu_0,0)$,
 $$\langle\tilde\cL_\mu u,u\rangle <  \langle -\Delta u,u\rangle,\quad\forall u\in C^\infty_0(\Omega)$$
 It follows by \cite[Theorem 10.2.2]{BS} that
  $$\lambda_{\mu,k}-\mu D_0^{-2}\leq \lambda_{0,k},\quad k=1,2,\cdots$$
 together with the upper bound  in \cite[(2.10)]{CgY}
 $$\lambda_{0,k}\leq \big(1+\frac4N\big)\lambda_{0,1}k^{\frac2N},$$
 we imply that
 $$\lambda_{\mu,k}\leq \big(1+\frac4N\big)\lambda_{0,1}k^{\frac2N}+\mu D_0^{-2},\quad k=1,2,\cdots.$$
 The proof is complete.\hfill$\Box$ \medskip

 \noindent{\bf Proof of Corollary \ref{cr 1.1}: }
   Recall that $(\lambda_{0,1},\phi_{0,1})$ is the first Dirichlet eigenvalue and associated eigenfunctions of Laplacian, then for $\mu>0$ and $N\geq 3$
 \begin{align*}
 \lambda_{\mu,1}\leq \int_\Omega |\nabla\phi_{0,1}|^2 dx+\mu\int_{\Omega} \frac{\phi_{0,1}^2}{|x|^2}dx\leq \lambda_{0,1} +\mu\|\phi_{0,1}\|_{L^\infty}^2\int_{\Omega} |x|^{-2} dx,
 \end{align*}
which, together with (\ref{e4}), implies that
  $$\lambda_{\mu,k}\leq \big(1+\frac4N\big)\Big(\lambda_{0,1} +\mu\|\phi_{0,1}\|_{L^\infty}^2\int_{\Omega} |x|^{-2} dx\Big)k^{\frac2N}.$$

For $\mu>0$ and $N\geq 2$, take test function $u=\phi_{0,1}\Gamma_\mu$ where $\Gamma_\mu(x)=|x|^{\tau_+(\mu)}\,{\rm with}\   \tau_+(\mu)=  -\frac{ N-2}2+\sqrt{\mu-\mu_0}$, which satisfying $\cL_\mu \Gamma_\mu=0\, {\rm in}\, \R^N\setminus\{0\}$. Then
 \begin{align*}
 \lambda_{\mu,1}&\leq
 \frac{\displaystyle\int_\Omega |\nabla u|^2 dx+\mu\displaystyle\int_{\Omega} \frac{u^2}{|x|^2}dx}{\|u\|^2_{L^2(\Omega)}}
 \\& = \lambda_{0,1} +\frac{\displaystyle\int_{\Omega} \nabla\phi_{0,1}^2\cdot \nabla\Gamma_\mu^2   dx }{2\displaystyle\int_{\Omega}\phi_{0,1}^2\Gamma_\mu^2 dx }
 \\&\leq \lambda_{0,1} +c_1^{-2}\tau_+(\mu)\|\phi_{0,1}\|^2_{C^1} \frac{\displaystyle\int_{\Omega}  |x|^{2\tau_+(\mu)-2} dx }{\displaystyle\int_{\Omega}\rho^2(x)  |x|^{2\tau_+(\mu)}  dx },
 \end{align*}
 where $c_1>0$ such that
 $$\phi_{0,1}(x)\geq c_1\rho(x),\quad\forall\, x\in\Omega.$$
 Together with (\ref{e4}), we obtain that
  $$\lambda_{\mu,k}\leq \big(1+\frac4N\big)\Big(\lambda_{0,1} +c_1^{-2}\tau_+(\mu)\|\phi_{0,1}\|^2_{C^1}  \frac{\displaystyle\int_{\Omega}  |x|^{2\tau_+(\mu)-2} dx }{\displaystyle\int_{\Omega}\rho^2(x)  |x|^{2\tau_+(\mu)}  dx }
  \Big)k^{\frac2N}.$$
 The proof is complete.\hfill$\Box$

 \section{Limit of Eigenvalues}

 Our proof of Weyl's limit of eigenvalues relies on estimates of the partition function
 $$Z_\mu(t):=\sum^{\infty}_{i=1} e^{-\lambda_{\mu,i}t},\quad t>0,$$
 which  can be written as
\begin{equation}\label{5.1-0}
Z_\mu(t)=\int_0^\infty e^{-\beta t} d\cN_\mu(\beta),
\end{equation}
 where $\cN_\mu(\beta)=\sum_{\lambda_{\mu,i}\leq \beta} 1$ is the usual counting function.
 Additionally, it has another standard formula for the partition function
\begin{equation}\label{5.1}
Z_\mu(t)=\int_{\Omega} p_{\mu,\Omega}(x,x,t) dx,
\end{equation}
where $p_{\mu,\Omega}$ is the heat kernel of Hardy operators $\cL_\mu$ in domain $\Omega\times \Omega\times(0,+\infty)$. More properties of heat kernel with Hardy potentials could be found in \cite{FMT,MT2}.   In the whole space $\Omega=\R^N$, we denote
$p_{\mu}=p_{\mu,\R^N}$. Particularly,
$$p_0(x,y,t)=p_0(x-y,t)=\frac1{(4\pi t)^{\frac N2}} e^{-\frac{|x-y|^2}{4t}}\ \ {\rm for}\ \  (x,y,t)\in \R^N\times \R^N\times(0,+\infty)$$
is the usual heat kernel of the Laplacian in $\R^N$.

The following lemma concerns the estimate of heat kernel, which is the essential part for the Weyl's limit.

\begin{lemma}\label{lm 4.1}
Let  $\mu\geq \mu_0$ and
$$\gamma_\mu(x,y,t)=p_0(x-y,t)-p_{\mu,\Omega}(x,y,t)\ \ {\rm for}\ \  (x,y,t)\in \R^N\times \R^N\times(0,+\infty). $$
Then\\
$(i)$ For  $\mu>0$, there exists $c_2>0$ independent of $\mu$ such that
$$0\leq \gamma_\mu(x,x,t)\leq  c_2\mu  |x|^{-2} t^{-\frac N2+1},\quad \forall\, x\in\Omega\setminus\{0\},\  \forall\,  t\in(0,1]. $$

\noindent $(ii)$ For  $\mu\in(\mu_0,0)$,   there exists $c_3>0$ such that
$$|\gamma_\mu(x,x,t) | \leq  c_3 |x|^{2\tau_+(\mu)-2}t^{-\frac N2+1},\quad \forall\, x\in\Omega\setminus\{0\},\  \forall\,  t\in(0,1].$$

\noindent $(iii)$ For  $\mu=\mu_0$,   there exists $c_4>0$ such that
$$|\gamma_\mu(x,x,t) | \leq  c_4  \Big( |x|^{-N+\frac{N+2}{8(N-2)}}t^{-\frac{N}{2}+1} +|x|^{-N+\frac18} t^{-\frac{N}{2}+\frac78}\Big),\quad \forall\, x\in\Omega\setminus\{0\},\  \forall\,  t\in(0,1].$$
\end{lemma}

\noindent{\bf Proof. }
{\it Case of $\mu>0$. } Note that elementary properties  infer directly that
$$p_{\mu,\Omega}(x,y,t)\leq p_{0,\Omega}(x,y,t) \leq p_{0}(x-y,t)  \quad {\rm in}\quad \Omega\times \Omega\times (0,+\infty)$$
and for fixed $y\in\Omega$, we have that
\begin{equation}\label{eq 4.1}
\left\{ \arraycolsep=1pt
\begin{array}{lll}
 \partial_t\gamma_\mu-\Delta \gamma_\mu=\frac{\mu}{|x|^2} p_\mu\quad \  &{\rm in}\ \   \Omega \times (0,+\infty),\\[2mm]
 \phantom{ \partial_t\gamma_\mu-\Delta    }
  \gamma_\mu=0\quad \ &{\rm{on}}  \ \, \partial\Omega \times (0,+\infty),\\[2mm]
 \phantom{ \partial_t\gamma_\mu-\Delta    }
 \gamma_\mu=0\quad&{\rm in}\ \ \Omega\times  \{0\}.
\end{array}
\right.
\end{equation}
 Note that $\gamma_\mu$ could be expressed by heat kernel, i.e. for  $y\in\Omega\setminus\{0\}$,
 \begin{align*}
 \gamma_\mu(x,y,t)  &=\int^t_0\int_{\Omega} p_{0,\Omega}(x,z,t-s) \frac{\mu   }{|z|^2} p_{\mu,\Omega}(z,y,s)dzds
 \\&\leq \mu \int^t_0\int_{\R^N} p_{0}(x-z,t-s) \frac{1 }{|z|^2} p_0(z-y,s)dzds
 \\&\leq\frac\mu{(4\pi)^{ N}} \int^t_0\frac1{(t-s)^{\frac N2} s^{\frac N2}} \int_{\R^N}  \frac{1 }{|z|^2} e^{-\frac{|x-z|^2}{4(t-s)}} e^{-\frac{|z-y|^2}{4s}} dzds.
 \end{align*}
 For   $N\geq 3$ and $x\in \Omega\setminus \{0\}$, we have that
 \begin{align*}
 \gamma_\mu(x,x,t)  & \leq\frac\mu{(4\pi)^{ N}} \Big(\int^t_0\frac1{(t-s)^{\frac N2} s^{\frac N2}}  \int_{\R^N\setminus B_{\frac{|x|}{2}}(x)} \frac{1 }{|z|^2} e^{-\frac{|x-z|^2}{4(t-s)}} e^{-\frac{|z-y|^2}{4s}} dzds
 \\  &\qquad\qquad\ \  + \int^t_0\frac1{(t-s)^{\frac N2} s^{\frac N2}}  \int_{ B_{\frac{|x|}{2}}(x)} \frac{1 }{|z|^2} e^{-\frac{|x-z|^2}{4(t-s)}} e^{-\frac{|x-z|^2}{4s}} dzds\Big)
 \\ & < \frac\mu{(4\pi)^{ N}} \Big(4|x|^{-2}\int^t_0\frac1{(t-s)^{\frac N2} s^{\frac N2}}  \int_{\R^N }   e^{-\frac14(\frac1{t-s}+\frac1{s}) |z|^2 }  dzds
 \\  &\qquad\qquad\ \  + \int^t_0\frac1{(t-s)^{\frac N2} s^{\frac N2}} e^{-\frac{ |x|^2}{16} \frac{t}{(t-s)s}   }  \int_{ B_{\frac{|x|}{2}}(0)} \frac{1 }{|z|^2}   dzds\Big)
\\ & \leq  \frac\mu{(4\pi)^{ N}} \Big( 2^{\frac N2+2} |x|^2t^{-\frac N2+1}  \int_{\R^N}e^{-|\tilde z|^2}d\tilde z
 \\  &\qquad\qquad\ \  +c_5 \frac{\omega_{_{N-1}}}{N-2} (\frac{|x|}2)^{N-2} \int^t_0\frac1{(t-s)^{\frac N2} s^{\frac N2}}  (\frac{ |x|^{-2}}{16} \frac{t}{(t-s)s} )^{-\frac N2}  ds   \Big)
\\ & =  \frac\mu{(4\pi)^{ N}}  |x|^{-2} t^{-\frac N2+1} \Big( 2^{\frac N2+2}  \int_{\R^N}e^{-|z|^2}dz +c_5\frac{\omega_{_{N-1}}}{N-2} 2^{N+2}     \Big),
 \end{align*}
 where $c_5>0$ such that
 $$e^{-a}\leq c_5 a^{-\frac N2}\quad{\rm for\ all}\ \  a\geq 1. $$

{\it Case of $\mu\in(\mu_0,0)$. }
It is known that for  $y\in\Omega\setminus\{0\}$,
 \begin{align*}
 \gamma_\mu(x,y,t)  &=\int^t_0\int_{\Omega} p_{0,\Omega}(x,z,t-s) \frac{\mu   }{|z|^2} p_{\mu,\Omega}(z,y,s)dzds.
 \end{align*}
 It follows from  \cite[Theorem 3.10]{MT2} that the corresponding heat kernel verifies that for $t\in(0,1]$
\begin{eqnarray*}
0\leq p_{\mu,\Omega}(x,y,t)  \le p_{\mu,\R^N}(x,y,t) \leq   c_6  (|x||y|)^{ \tau_+(\mu) } t^{-\frac N2} e^{- c_7\frac{|x-y|^2}{t} },
\end{eqnarray*}
where $c_6, c_7>0$ and $\tau_+(\mu)<0$. We see that there is no order for $p_{\mu,\Omega}$ and $p_0$.
Thus,  for $x\in\Omega\setminus\{0\}$ and $t\in(0,1]$
\begin{align*}
| \gamma_\mu(x,x,t) | &\leq  |\mu|\int^t_0\int_{\Omega} p_{0,\Omega}(x,z,t-s) \frac{1}{|z|^2} p_{\mu,\Omega}(z,x,s)dzds
 \\&\leq \frac{c_6   |\mu| }{(4\pi)^{ \frac{N}2 }} |x|^{\tau_+(\mu)} \int^t_0\frac{1}{(t-s)^{\frac N2}s^{\frac N2}  } \int_{\R^N}   |z|^{ \tau_+(\mu)-2 }   e^{-\frac{|x-z|^2}{4(t-s)}-\frac{c_7|x-z|^2}{ s} } dzds
\\ & \leq  \frac{c_6   |\mu| }{(4\pi)^{ \frac{N}2}} |x|^{\tau_+(\mu)}  \Big(\int^t_0\frac{1}{(t-s)^{\frac N2} s^{\frac N2}  } \int_{\R^N\setminus B_{\frac{|x|}{2}}(x)}  |z|^{ \tau_+(\mu)-2 }   e^{-\frac{|x-z|^2}{4(t-s)} -\frac{c_7|x-z|^2}{ s} } dzds
 \\  & \qquad\qquad\qquad  \quad + \int^t_0\frac{1}{(t-s)^{\frac N2} s^{\frac N2} }  \int_{ B_{\frac{|x|}{2}}(x)}   |z|^{ \tau_+(\mu)-2 }   e^{-\frac{|x-z|^2}{4(t-s)}-\frac{c_7|x-z|^2}{ s}  } dzds\Big)
 \\ & <  \frac{c   |\mu| }{(4\pi)^{ \frac{N}2}} |x|^{\tau_+(\mu)} \Big((\frac{|x|}2)^{\tau_+(\mu)-2}\int^t_0\frac{1}{(t-s)^{\frac N2} s^{\frac N2} }  \int_{\R^N }   e^{- \frac{|x-z|^2}{4(t-s)}-\frac{c_7|x-z|^2}{ s}    }  dzds
 \\  & \qquad\qquad\qquad  \qquad  + \int^t_0\frac{1}{(t-s)^{\frac N2} s^{\frac N2} }  e^{- \frac{|x|^2}{16(t-s)} -\frac{c|x|^2}{ s}   }  \int_{ B_{\frac{|x|}{2}}} |z|^{\tau_+(\mu)-2}  dzds\Big)
\\ & \leq   \ \frac{c_6   |\mu| }{(4\pi)^{ \frac{N}2}} |x|^{\tau_+(\mu)} \Big(  (\frac{|x|}2)^{\tau_+(\mu)-2}  \int_{\R^N}e^{-|z|^2}dz \int_0^t (c_7t+(4-c_7)s)^{-\frac{N}{2}} ds
 \\  &  \qquad\qquad\  +\frac{\omega_{_{N-1}}}{N-2} (\frac{|x|}2)^{N+\tau_+(\mu)-2}c_5  \int^t_0\frac1{(t-s)^{\frac N2} s^{\frac N2} }  \Big(\frac{ |x|^{2}}{4(t-s)}+c_7\frac{ |x|^{2}}{s} \Big)^{-\frac N2}  ds   \Big)
\\ & \leq   \frac{c_8   |\mu| }{(4\pi)^{ \frac{N}2}}   \Big( 2^{2-\tau_+(\mu)}   \int_{\R^N}e^{-|z|^2}dz +c_9 \frac{\omega_{_{N-1}}}{N-2} 2^{N-\tau_+(\mu)+2} c^{-\frac N2}   \Big)|x|^{2\tau_+(\mu)-2}t^{-\frac{N}{2}+1} ,
\end{align*}
 where $2\tau_+(\mu)-2>-N$ for $\mu\in(\mu_0,\,0)$ and $c_8,c_9>0$.\smallskip

 {\it Case of $\mu=\mu_0$. }  When $\mu=\mu_0$, the above inequality holds true, but
 the factor $|x|^{2\tau_+(\mu_0)-2}=|x|^{-N}$ is non-integrable in $\Omega$. So we have to modify
 the above estimates.

   From  \cite[Theorem 1.1]{FMT} we have that
 $t\in(0,1]$
\begin{eqnarray*}
0\leq p_{\mu,\Omega}(x,y,t)  \le    c_{10}  (|x||y|)^{ \frac{2-N}{2} } t^{-\frac N2} e^{- c_{11}\frac{|x-y|^2}{t} }.
\end{eqnarray*}
 We choose
 $$r=\frac14 |x|^{1-\frac{1}{4(N-2)}}\quad{\rm and}\quad \theta=\frac14,$$
then
 we have that
 \begin{align*}
| \gamma_\mu(x,x,t) | &\leq    \frac{c_{10}   |\mu_0| }{(4\pi)^{ \frac{N}2}} |x|^{ \frac{2-N}{2}}  \Big(\int^t_0\frac{1}{(t-s)^{\frac N2} s^{\frac N2}  } \int_{\R^N\setminus B_r(x)}  |z|^{  -\frac{N+2}{2} }   e^{-\frac{|x-z|^2}{4(t-s)} -\frac{c_{11}|x-z|^2}{ s} } dzds
 \\  & \qquad\qquad\qquad  \quad + \int^t_0\frac{1}{(t-s)^{\frac N2} s^{\frac N2} }  \int_{ B_r(x)}   |z|^{ -\frac{N+2}{2} }   e^{-\frac{|x-z|^2}{4(t-s)}-\frac{c_{11}|x-z|^2}{ s}  } dzds\Big)
 \\ & <  \frac{c_{10}   |\mu_0| }{(4\pi)^{ \frac{N}2 }} |x|^{\frac{2-N}{2}} \Big(r^{-\frac{N+2}{2}}\int^t_0\frac{1}{(t-s)^{\frac N2} s^{\frac N2} }  \int_{\R^N }   e^{- \frac{|x-z|^2}{4(t-s)}-\frac{c_{11}|x-z|^2}{ s}    }  dzds
 \\  & \qquad\qquad\qquad  \qquad  + \int^t_0\frac{1}{(t-s)^{\frac N2} s^{\frac N2} }  e^{- \frac{|x|^2}{16(t-s)} -\frac{c_{11}|x|^2}{ s}   }  \int_{ B_r(0)} |z|^{-\frac{N+2}{2}}  dzds\Big)
\\ & \leq   \ \frac{c_{10}   |\mu_0| }{(4\pi)^{ \frac{N}2 }} |x|^{\frac{2-N}{2}} \Big( r^{-\frac{N+2}{2}}  \int_{\R^N}e^{-|z|^2}dz \int_0^t (c_{11}t+(4-c_{11})s)^{-\frac{N}{2}} ds
 \\  &  \qquad\qquad\  +\frac{\omega_{_{N-1}}}{N-2} r^{\frac{N-2}{2}}c_{14}  \int^t_0\frac1{(t-s)^{\frac N2} s^{\frac N2} }  \Big(\frac{ |x|^{2}}{4(t-s)}+c_{11}\frac{ |x|^{2}}{s} \Big)^{-\frac {N-\theta}2 }  ds   \Big)
\\ & \leq    c_{12}   |x|^{\frac{2-N}{2}}r^{-\frac{N+2}{2}} t^{-\frac{N}{2}+1} + c_{13}|x|^{1+\theta-\frac{3}{2}N}r^{\frac{N-2}{2}}  t^{-\frac{N-\theta}{2}}\int_0^t(t-s)^{-\frac\theta2}s^{-\frac\theta2} ds
\\ & =  c_{12}   |x|^{-N+\frac{N+2}{8(N-2)}}t^{-\frac{N}{2}+1} +c_{13}|x|^{-N+\frac18} t^{-\frac{N}{2}+\frac78} ,
\end{align*}
where
$c_{12},c_{13}>0$ and $c_{14}>0$ verifies that
 $$e^{-a}\leq c_{14} a^{-\frac {N-\frac14}2 }\quad{\rm for\ all}\ \  a\geq 1. $$
 This  completes the proof. \hfill$\Box$ \medskip

\noindent{\bf Proof of Theorem \ref{teo limit}. }
  From Lemma \ref{lm 4.1},  we see that for $t\in(0,1)$ and $\mu\geq\mu_0$,
\begin{eqnarray*}
t^{\frac N2}\Big| \int_{\Omega}
 p_{\mu,\Omega}(x,x,t)dx -\int_{\Omega}
p_0(0,t)dx\Big|&\leq  &  t^{\frac N2} \int_{\Omega} |\gamma_{\mu}(x,x,t)| d x\\&\leq&  c_2 t ^{\frac78}    \int_{\Omega} \max\{  |x|^{ -2},|x|^{2\tau_+(\mu)-2}\}dx
\\ &\to& 0\ \ \ \text{ as $t \rightarrow 0^+$},
\end{eqnarray*}
then we obtain that
$$\lim_{t\to0^+} t^{\frac N2} |Z_\mu(t)-Z_0(t)|=0,$$
where $Z_0$ is the partition function for Laplacian and direct computation shows that
$$
\lim_{t\to0^+}t^{\frac N2} Z_0(t)=\frac{|\Omega| }{(4\pi )^{\frac N2}}.
$$
So there holds
\begin{equation}\label{5.2}
\lim_{t\to0^+} t^{\frac N2} Z_\mu(t)= \frac{|\Omega| }{(4\pi )^{\frac N2}}.
\end{equation}
From (\ref{5.2}) and Karamata's Tauberian theorem \cite[Theorem 10.3]{Sim},   we get that
$$\cN_\mu(\beta) \beta^{-\frac N2}\longrightarrow \frac{|\Omega| }{\Gamma(\frac N2+1)(4\pi )^{\frac N2}}\quad{\rm as}\ \beta\to+\infty,$$
which, taking $\beta=\lambda_{\mu,k}$ and $\cN(\lambda_{\mu,k})=k$, implies that
\begin{equation}\label{www-1}
\lim_{k\to+\infty}\lambda_{\mu,k}^{-\frac{N}2}k=\frac{|\Omega| }{\Gamma(\frac N2)\frac N2 (4\pi )^{\frac N2}} =\frac{|\Omega| |B_1|}{  (2\pi )^{N}},
\end{equation}
where $$|B_1|=\frac1N\omega_{_{N-1}}=\frac{2\pi^{\frac N2}}{N\Gamma(\frac N2)}.$$
Note that (\ref{www-1}) is equivalent to   (\ref{lim 1}) and the proof is completed.  \hfill$\Box$\medskip

\bigskip

\medskip
 \noindent{\small {\bf Acknowledgements:}  H. Chen is is supported by NNSF of China, No: 12071189 and 12001252,
by the Jiangxi Provincial Natural Science Foundation, No: 20202BAB201005, 20202ACBL201001. F. Zhou is supported by Science and
Technology Commission of Shanghai Municipality (STCSM), Grant No. 18dz2271000 and also supported by
NSFC (No. 11431005).

\end{document}